\documentclass[12pt]{iopart}

\usepackage{amsthm}
\usepackage{amsfonts}

\begin{document}

\title{Local Invertibility of Integral Operators with Analytic Kernels}

\author{Nikolay Balov}

\address{
				Department of Biostatistics and Computational Biology \\
				University of Rochester Medical Center \\
				601 Elmwood Ave \\
				Rochester, NY-14642 \\
		                \ead{nikolay\_balov@urmc.rochester.edu}
}

\maketitle

\begin{abstract}
The invertibility of integral linear operators is a major problem of both theoretical and practical importance. 
In this paper we investigate the relation between an operator invertibility and the rank of its integral kernel to develop a local inverse theory. An operator is called locally invertible provided that any function can be recovered from its transformed image if the latter is known in an arbitrary open subset of its domain, i.e., if its image is known locally. It turns out that the local invertibility holds for any analytical kernel whose Taylor functions are linearly independent in any open subset of their domain - the so-called local linear independence condition.  
We also establish an equivalence between local linear independence and the so-called full rank a.e. property. The latter can be described as follows: for any finite, random sample of points, the square matrix obtained by applying, pairwise, the kernel function on them, has full rank almost surely. As an illustration, we show that the geodesic distance function on a sphere in more than one dimensions is of full rank a.e., in contrast to the Euclidean distance which is not.
\end{abstract}



\section{Problem Formulation and Motivation}

A large class of linear inverse problems arises from the Fredholm integral equation of the first kind 
\begin{equation}\label{eq:ram_fredholm}
\int_V \psi(x,y) f(y)dy = g(x), \textrm{ }x\in U,
\end{equation}
where $U$ and $V$ are open sets in $\mathbb{R}^n$ and $f:V\to\mathbb{R}$, $g:U\to\mathbb{R}$, 
and $\psi:U\times V\to\mathbb{R}$ are some functions. 
The so-called kernel $\psi$ defines a linear operator $A_{\psi}$, allowing (\ref{eq:ram_fredholm}) to be conveniently written as $A_{\psi}f = g$. Typically, $f$ and $g$ are elements of the $L^2$-spaces $L^2(V)$ and $L^2(U)$, while $\psi$ is a Hilbert-Schmidt kernel, in which case $A_{\psi}$ is a compact linear operator. Notable examples of $A_{\psi}$ are the Laplace transformation in $\mathbb{R}$ with kernel $\psi(x,y) = \exp(-xy)$ and the closely related Fourier transformation, $\psi(x,y) = \exp(- i xy)$. Without being specific, we shall assume that (\ref{eq:ram_fredholm}) is formulated in $L^p$ sense for some $p\ge 1$. 

Let us assume that for given $g$, (\ref{eq:ram_fredholm}) is solvable. Under some favorable conditions, a numerical solution of (\ref{eq:ram_fredholm}) can be obtained using a quadrature method by sampling the functions $f$, $g$ and $\psi$ at discrete sets of points. Let $\bar U=\cup_{i=1}^k \bar U_i$ and $\bar V=\cup_{j=1}^k \bar V_j$ be partitions of the closures $\bar U$ and $\bar V$ such that $U_i$ and $V_j$ are open and disjoint. 
Let also $\{x_i\}_{i=1}^k$ and $\{y_i\}_{i=1}^k$ be two sets of points such that $x_i\in U_i$ and $y_j\in V_j$. 
Then, we can discretize equation (\ref{eq:ram_fredholm}) to 
\begin{equation}\label{eq:ram_fredholm_discrete}
\sum_{j=1}^k \psi(x_i, y_j) f(y_j) vol(V_j) = g(x_i), \textrm{ }i=1,...,k, 
\end{equation}
a system of linear equations. 
System (\ref{eq:ram_fredholm_discrete}) can be solved for $f_k=f(y_j)_{j=1}^k$, only if the matrix $A_k = \{\psi(x_i,y_j)\}_{i,j=1}^{k}$ is invertible, or equivalently, has full rank $k$; then the solution is $\hat{f}_k(y) = \sum_{j=1}^k \hat{f}_{k,i} 1_{V_j}(y)$, for the vector $(\hat{f}_{k,j})_{j=1}^k = W_k A_k^{-1} g_k$, where $W_k = diag(1/vol(V_j))$ and $g_k=(g(x_i))_{i=1}^k$. Let us assume that $A_{\psi}$ is invertible and for a sequence of nets $\{x_i,y_j\}_{i,j=1}^{k}$, $k=1,2,...$, $A_k$ converges uniformly to $A_{\psi}$ and $g_k$ converges uniformly to $g$. Then, any convergent subsequence of $\hat{f}_k$ will approximate the solution of (\ref{eq:ram_fredholm}) (more details are given in the concluding section). Even though often, such as for compact $A_{\psi}$, $A_{\psi}^{-1}$ is unbounded and convergent subsequence of $\hat{f}_k$ may not be available, the solvability of equation (\ref{eq:ram_fredholm_discrete}) is an interesting question by itself that is relevant to many problems with discrete nature. 

Usually, the sample points $(x_i,y_j)$ are chosen in a deterministic way, for example, according to uniform spacing. However, in order to formulate a more general problem regarding the invertibility of $A_k$, we introduce the following stochastic setting. Let the points $x_i$ and $y_j$ be chosen independently by some continuous distributions 
on $U$ and $V$, respectively. Then, the rank of $A_k$ becomes a $2k$-dimensional continuous random variable. 
We are interested in the question: what property of the kernel $\psi$ guarantees that for every $k$, $A_k$ has full rank almost surely (a.s.)? If this is the case, we shall say that $\psi$ has full rank almost everywhere (a.e.) in its domain.

The main focus of this paper is the characterization of the full rank a.e. kernels. We first need to define them properly. 
A possible approach is to simply reverse the standard, finite rank definition and say that a kernel $\psi$ has full rank, if it does not admit a finite sum representation $\psi(x,y) = \sum_{l=1}^m \phi_l(x)\xi_l(y)$ for linearly independent $\phi_l$ - finite rank linear operators are considered in all textbooks on the subject, see for example $\$70$ in \cite{riesz-nagy}. This approach however will be inadequate, for thus defined full rank kernels may not be of full rank in ``a.e.'' sense, as the following example shows. In $L^2[0,1]$, consider the functions 
\begin{equation*}
  c_s(x) = 1_{[\frac{s}{s+1}, \frac{s+1}{s+2})}(x) \textrm{, } s \ge 0, 
\end{equation*}
and let $\psi(x,y) = \sum_{s=0}^{\infty} \frac{1}{s!} c_s(x) y^s$ - we call analytic in $y$ all kernels $\psi$ that admit such Taylor expansion. Note that, as defined, $c_s$ are linearly independent in $[0,1]$ and thus, $\psi(x,y)$ has no finite rank in sense of the standard definition. On the other hand, for any $k$ points $x_i$ and $y_j$, sampled uniformly in $[0,1]$, there is a positive probability for $A_k$ to have rank $1$, and hence, $\psi$ is not of full rank a.e. in $[0,1]$. 
Evidently, for the latter to hold, $c_s$ need to be linearly independent in a neighborhood of any point of their domain. We shall call the latter requirement on $c_s$ local linear independence. 
In the example above, $c_s$ are linearly independent but not locally linear independent. 
The main result in the paper is about rigorously showing the equivalence between the full rank a.e. property of analytical kernels and the local linear independence of their Taylor functions. 

Local linear independence is a necessary condition for the solvability of some classes of integral equation problems. Our interest in this topic started with the study of the so-called variance linear operator in the domain of probability density functions. More generally, in (\ref{eq:ram_fredholm}), let $U$ and $V$ be open subsets of $\mathbb{R}^n$, and $\psi(x,y)=d(x,y)$ be a distance function in $\mathbb{R}^n$. Consider the operator $C_U$
$$
(C_{U}f)(x) = \int_V d^2(x,y)f(y)dy \textrm{, } x\in U, 
$$
where $f\in \mathbb{P}(V)$, a class of measurable functions in $V$ with finite moments of all orders. If $d^2$ is analytic in $y\in V$, $d^2(x,y) = \sum_{s=0}^{\infty} c_s(x) y^s$, then the local linear independence of $c_s$ in $U$ will guarantee the invertibility of $C_{W}$ for any open $W\subset U$ (for more details see the concluding section). Consequently, there will be a one-to-one map between $\mathbb{P}(V)$ and the space of variance functions from the image $C_{W}(\mathbb{P}(V))$. With other words, any element of $\mathbb{P}(V)$ can be recovered if we know its variance function in an arbitrary open domain. This is what we tentatively call local invertibility. On the other hand, if $c_s$ are just linearly independent in $U$, then one needs to know the variance of $f$ everywhere in $U$ in order to recover it. 
Note that if only finitely many $c_s$ are non-zero, then the system $\{c_s\}_s$ can be neither linearly nor locally linear independent and no recovery is possible; this is the case with the Euclidean distance, $d^2(x,y) = ||x-y||^2$. As we show later in this paper, an example of full rank a.e. kernel is the squared spherical distance in $\mathbb{R}^n$, $d^2(x,y) = arccos^2(\frac{x.y}{||x||||y||})$. At this point we shall speculate that by establishing a one-to-one correspondence between $\mathbb{P}(V)$ and $C_{U}(\mathbb{P}(V))$, one obtains new potential tools for studying probability distributions, especially such on non-Euclidean metric spaces. 

The local invertibility of integral operators may have direct practical implications. This is evident from the fact that 
the integral equation (\ref{eq:ram_fredholm}) and its discrete analog (\ref{eq:ram_fredholm_discrete}) provide an abstract  formulation for many remote sensing inverse problems, where the goal is to study some directly unobserved phenomena $f$ from its integral characteristic $g=A_{\psi}f$, which is measurable. Then, at least in principle, the local invertibility of $A_{\psi}$ allows the measurements on $g$ to be taken in an arbitrary domain of the measurement space and still reconstruct $f$. It is no coincidence then that all important in practice integral transformations, such as the Laplace transformation, are locally invertible. 

The full rank a.e. property of integral kernels may also be of importance to some families of discrete inverse problems. 
Often, linear systems of equations of type (\ref{eq:ram_fredholm_discrete}), with continuously changing sample points, arise naturally in problems such as the following variant of the $n$-body problem in astronomy. Let us assume that $\psi$ represents a gravitation field, $f$ is a mass distribution and $g$ is an observed effect of the field generated by the mass, for example, acceleration. If the gravitation field is generated by $k$ otherwise freely moving objects with positions $\{x_i\}_{i=1}^k$, one may be interested in reconstructing their mass from an instant observed effect of the accumulated force on each of them, which yields a type (\ref{eq:ram_fredholm_discrete}) problem ($y_i$ are taken to be $x_i$). In this setting, the points $\{x_i\}_{i=1}^k$ will behave as if coming from a continuous probability distribution on the product space and the fact that the matrix $A_k$ has full rank a.s. and thus, (\ref{eq:ram_fredholm_discrete}) is solvable, has an immediate significance. 

The main body of the paper is organized in three sections. We begin, Section \ref{sec:finite_rank}, with a systematization of the finite rank condition for analytic kernels by presenting some, not necessarily new, equivalent conditions for the latter (Proposition \ref{theorem:rank_equivals}). Then, in Section \ref{sec:lli}, we introduce the local linear independence and derive an important necessary condition (Lemma \ref{lemma:lli_cond}). In Section \ref{sec:full_rank}, we define the full rank a.e. notion and present some necessary and sufficient conditions for it that hold for analytic kernels (Proposition \ref{theorem:psisum_fullrank} and Corollary \ref{cor:analytic_full_rank_cond}). Also there, presented is a more involved example of a full rank a.e. kernel - the spherical distance in two or more dimensions (Proposition \ref{theorem:full_rank_sphere}). Finally, we conclude with a discussion on the usefulness and limitations of the presented results, and their possible further development.

\section{Analytic Kernels with Finite Rank}\label{sec:finite_rank}

In this section we present some necessary and sufficient conditions for an analytic kernel to have finite rank. 
These conditions are based upon the fundamental relation between kernel rank and functional linear independence. 
Our goal here is not to present new facts, but to systematize some basic results and present them in a form needed in the following sections. 

Let $U$ be an open subset of $\mathbb{R}^n$. 
A collection of functions $f_s:U \to \mathbb{R}$ is said to be linear independent in $U$ if 
$\sum_s \alpha_s f_s(x) = 0$, for almost all (in Lebesgue measure sense) $x\in U$ only if $\alpha_s = 0$ for all $s$.
It is easy to check the following necessary condition for linear independence. 
\newtheorem{lemma_linindep}{Lemma}
\begin{lemma_linindep}
If the functions $f_1(x)$, ...,$f_k(x)$, $k\ge1$, are linearly independent in $U$, then there exist 
$x_1,...,x_k$ in $U$, such that $rank(\{f_i(x_j)\}_{i,j=1}^{k}) = k$.
\label{lemma:linindep_functions}
\end{lemma_linindep}

In the theory of linear operators, finite rank are said to be all kernels admitting representations in the form $\sum_{j=1}^k f_j(x)g_j(y)$, for linearly independent $f_j$ and $g_j$. We shall adopt, however, a different definition, based on the notion of matrix rank, which is equivalent to the former. The reason of not using the standard one is that our finite rank definition can be naturally extended to characterize full rank kernels as well and thus, provides more intuitive and consistent approach. 
\newtheorem{ram_def1}{Definition}
\begin{ram_def1}
We say that a function $\psi:U\times V \to \mathbb{R}$, $U,V\subset\mathbb{R}^n$, 
has rank $k$ and write $rank(\psi) = k$ if 
for any $m\in\mathbb{N}$, $x_i\in U$ and $y_j \in V$, i, j=1,...,m, 
$$
rank( \{\psi(x_i, y_j)\}_{i,j=1}^{m} ) \le k, 
$$
and $k$ is the smallest number with this property. 
\label{ram:rank_function}
\end{ram_def1}

Recall that a function $f\in C^{\infty}(U)$ is real analytic in $U$ if it admits power series expansion in a neighborhood of any point $x\in U$. Next we extend this definition to bivariate functions or kernels in the linear operator theory context. 
Hereafter, we assume that $U$ and $V$ are open connected subsets of $\mathbb{R}^n$, although more general settings may be possible. 
We say that $\psi:U\times V \to \mathbb{R}$ is analytic in $V$, if about any $p\in V$ one can write 
$$
\psi(x,y) = \sum_{m=0}^{\infty} \sum_{s:\sum_{i=1}^n s_i = m} c_{s_1...s_n}(x) (y_1-p_1)^{s_1}...(y_n-p_n)^{s_n}
\textrm{, } 
$$
for some functions $c_{s_1...s_n}$ in $U$, such that the series on the right converges to $\psi(x,y)$ for any $(x,y)\in U\times V$. 
For the sake of brevity, we shall write the above expansion as 
\begin{equation}\label{eq:psi_analytic}
\psi(x,y) = \sum_{s}^{\infty} c_s(x) (y-p)^s := \sum_{l=1}^{\infty} \sum_{s:[s]=l} c_s(x) (y-p)^s
\end{equation}
where $s=(s_1,...,s_n)$ is multi-index and $[s]:= s_1+...+s_n$. Any analytic $\psi(x,y)$, therefore, is infinitely differentiable in $y$ about any $p\in V$, and the series (\ref{eq:psi_analytic}) is in fact the Taylor series of $\psi$ in $y$ at $p$, thus giving $c_s(x) = \frac{1}{s!}\frac{\partial^s \psi}{\partial y^s}|_p$. We also silently assume that for all $x\in U$ and $p\in V$, the closed ball with center $p$ and radius the convergence radius of the series (\ref{eq:psi_analytic}), encloses the whole $V$. With other words, we assume that for any $p\in V$, $\psi$ admits an analytic expansion (\ref{eq:psi_analytic}) almost everywhere in $V$, except, eventually, a set of Lebesgue measure zero. If in (\ref{eq:psi_analytic}), $c_s(x)$ are analytic in $x\in U$ as well, then $\psi$ is said to be analytic in both arguments. 
The main results in this paper establish connections between the rank of $\psi$ and some conditions on its Taylor functions $c_s$. 

We say that the vector space $span\{c_1, c_2,...\}$, that is, the space of all finite linear combinations of $c_s$, has a finite basis of size $m$, if there exist functions $c_{s_1}$, ..., $c_{c_m}$, such that all $c_s$'s are their linear combinations, i.e. $c_s\in span\{c_{s_1},...,c_{s_m}\}$ for all $s$.
The next result gives some necessary and sufficient conditions for $\psi$ to have finite rank.
\newtheorem{rank_equivals}{Proposition}
\begin{rank_equivals}
For any function $\psi : U\times V \to \mathbb{R}$ that is analytic in $V$, 
the following three conditions are equivalent
\begin{enumerate}
\item[(1)] $rank(\psi) = k$.
\item[(2)] $\psi(x,y) = \sum_{j=1}^k \phi_j(x)\xi_j(y)$ for linearly independent 
functions $\phi_j:U\to\mathbb{R}$ and $\xi_j\in C^{\infty}(V)$, j=1,...,k.
\item[(3)] Let for arbitrary $p\in V$, $c_s$ are the Taylor functions of $\psi$ from Eq. (\ref{eq:psi_analytic}). Then, $span\{c_1, c_2,...\}$ has a finite basis of size $k$ in $U$.
\end{enumerate}
\label{theorem:rank_equivals}
\end{rank_equivals}

Note that without the analytical requirement on $\psi$ only the link $(2) \Rightarrow (1)$ seems to be evident, because the connection $(1) \Rightarrow (2)$ goes through the Taylor expansion of $\psi$. Also, if the condition $(3)$ holds for one $p\in V$ then it inevitably holds for any other point in $V$, and therefore, to show that $\psi$ has finite rank, one only needs to show that its Taylor functions about one particular point have a finite basis. 

If $\psi$ is symmetric and the representation in condition (2) of the proposition holds, then $\xi_i$ are linear combinations of $\phi_i$.Indeed, if we choose $x_j$, $j=1,...,k$, such that $rank(\{f_i(x_j)\}_{i,j=1}^{k}) = k$, 
then, for any $y$ we can solve the system 
$$
\sum_{i=1}^k \phi_i(x_j)\xi_i(y) = \sum_{i=1}^k \phi_i(y)\xi_i(x_j), j=1,...,k
$$
to obtain $\xi_i(y) = \sum_{i=1}^k \gamma_{ij} \phi_j(y)$. 

Another immediate corollary of Proposition \ref{theorem:rank_equivals} is that for analytical $\psi$ with finite rank, $rank(\psi) = k$, the integral operator 
$$
A_{\psi}: f \mapsto \int_V \psi(x,y)f(y)dy, \textrm{ } f:U\to \mathbb{R}, 
$$
is of finite rank in sense of the standard definition from the integral operator theory. 

Next, we present probably the simplest application of Proposition \ref{theorem:rank_equivals} showing the finite rank of the Euclidean metric.
Consider the Euclidean square distance in $\mathbb{R}^n$, 
$
\psi(x,y) = ||x-y||^2 = \sum_{s=1}^n (x_{[s]}-y_{[s]})^2,
$
where $x_{[s]}$ are the components of $x\in \mathbb{R}^n$.
Because of the global representation 
$$
\psi(x,y) = \sum_{s=1}^n x_{[s]}^2 * 1 - 2x_{[1]} y_{[1]} ... - 2 x_{[n]} y_{[n]} + 1 * \sum_{s=1}^n y_{[s]}^2,
$$
by Proposition \ref{theorem:rank_equivals}, $rank(\psi) = n+2$, i.e. 
the Euclidean metric has finite rank of two more than the number of dimensions.

Another obvious example of a finite rank kernel is the circular square distance, $\psi(x,y) = (x-y)^2$, $x,y\in \mathbb{S}^1 \equiv [0,2\pi]$. At the end of this paper we shall show that the spherical square distance in more than one dimensions is never of finite rank. 

\section{Local Linear Independence}\label{sec:lli}

In this section we introduce and describe the so-called local linear independence of a collection of functions, a property stronger than their linear independence. As it turns out, the rank of a bivariate analytic function is related to the local linear independence of its Taylor functions, which motivates the introduction of the latter notion. We start with its definition. 
\newtheorem{ram_def2}[ram_def1]{Definition}
\begin{ram_def2}
A collection of functions $f_s:U \to \mathbb{R}$ is said to be locally linear independent in $U$, if 
for any open subset $W$ of $U$ we have that 
$\sum_s \alpha_s f_s(x) = 0$ almost everywhere in $W$ if and only if $\alpha_s = 0$ for all $s$.
\label{def:local_linear_independence}
\end{ram_def2}
Apparently, local linear independence implies linear independence. 
The reverse however is not true - a collection of functions can be linear independent in one open region but dependent in another. 
This is exactly the case with the functions 
\begin{equation*}
  f_s(x) = \left\{
    \begin{array}{rl}
      e^{-s/x} & \textrm{if } x > 0,\\
      0 & \textrm{if } x \le 0,
    \end{array} \right.
\end{equation*}
which are not locally linear independent, because they are linearly dependent on the negative half of the real line.

Naturally, the most important example of local linear independence is presented by the power function. 
\newtheorem{lemma_lli_powers}[lemma_linindep]{Lemma}
\begin{lemma_lli_powers}\label{lemma:lli_powers}
For any distinct multi-indices $s_1,s_2,...$, the power functions $\{x^{s_l}\}_{l=1}^{\infty}$ are locally linear independent in $\mathbb{R}^n$.
\end{lemma_lli_powers}

In fact, the local linear independence of the power functions can be generalized to any collection of analytic functions - there is an equivalence between local linear independence and linear independence for any collection of analytic functions. The reason is that any analytic function vanishing in an open region, vanishes everywhere in its definition domain. 
Consequently, if a linear combination of real analytic functions vanishes in an open set, being an analytic function itself, it vanishes everywhere in its domain, and therefore the collection can not be linearly independent. 

Next, we state several elementary properties of local linear independence that we are going to use later. 
Note that any collection of translated power functions $\{(x-a_l)^{s_l}\}_{l=1}^{\infty}$ is also locally linear independent. In fact, any bijection that have continuous inverse preserves local linear independence.  
\newtheorem{lemma_lli_bijection}[lemma_linindep]{Lemma}
\begin{lemma_lli_bijection}\label{lemma:lli_bijection}
Let $U$ and $W$ be open subsets of $\mathbb{R}^n$ and $\mathbb{R}^k$, $k \le n$, respectfully, and $h:U\to W$ be a transformation with the property: the image of any open subset in $U$ contains an open subset in $W$. 
Then, for any collection $f_s: W \to \mathbb{R}$ of locally linear independent functions in $W$, $f_s\circ h$ are locally linear independent in $U$. 
\end{lemma_lli_bijection}
The assumption of Lemma \ref{lemma:lli_bijection} is satisfied by any $h:U\to W$ that has a full rank Jacobian in $U$, that is, since $k\le n$, $rank(J_h(p))=k$, $\forall p\in U$. 

Although, in general, local linear independence is not preserved by forming linear combinations, there are some special cases when this is true. 
\newtheorem{lemma_lli_lin_comb}[lemma_linindep]{Lemma}
\begin{lemma_lli_lin_comb}\label{lemma:lli_lin_comb}
If $f_s: W \to \mathbb{R}$, s=1,...,n, are locally linear independent functions and 
$$
g_1 \in span\{f_1\} \textrm{, } g_s \in span\{f_1,...,f_s\} \setminus span\{f_1,...,f_{s-1}\} \textrm{, } s > 0,
$$ 
then $g_s$ are also locally linear independent. 
\end{lemma_lli_lin_comb}

We continue with an important connection between the local linear independence of finite sets of functions and the rank of matrices sampled from these functions. The result is analogous to that in Lemma \ref{lemma:linindep_functions} and gives a necessary condition for local linear independence.
\newtheorem{lemma_lli_cond}[lemma_linindep]{Lemma}
\begin{lemma_lli_cond}\label{lemma:lli_cond}
If the functions $f_1(x)$, ...,$f_k(x)$ are locally linear independent in $U$ and 
$U_i$, $i=1,...,k$, are arbitrary (non-empty) open subsets of $U$, then there exist 
$x_1,...,x_k$, $x_i\in U_i$, such that {$rank(\{f_i(x_j)\}_{i,j=1}^{k}) = k$}.
\end{lemma_lli_cond}
The condition in Lemma \ref{lemma:lli_cond} is only necessary. Indeed, it is clearly not-sufficient for the function $f(x) = 1_{\{x\textrm{ is rational}\}}$, $x\in (0,1)$, to be linearly independent, less locally linear independent. However, if in addition functions $f_i$ are continuous, then the condition becomes sufficient as well necessary. 

\section{Full Rank Almost Everywhere Kernels}\label{sec:full_rank}

We continue with an introduction of the full rank notion. As it shall become clear shortly, the definition is a natural extension of that of the finite rank and is closely related to local linear independence. 
\newtheorem{ram_def3}[ram_def1]{Definition}
\begin{ram_def3}
A bivariate function $\psi:U\times V \to \mathbb{R}$ is said to have full rank almost everywhere ($a.e.$) in $U\times V$ if 
for any number $k\in\mathbb{N}$ and open sets $U_i\subset U$ and $V_j\subset V$, $i,j=1,...,k$, 
there exist $x_i\in U_i$ and $y_j\in V_j$ such that 
$$
rank( \{\psi(x_i, y_j)\}_{i,j=1}^{k} ) = k.
$$
\label{def:fullrank_function}
\end{ram_def3}
To say that $\psi$ has no finite rank according to Definition \ref{ram:rank_function}, is not equivalent to the full rank a.e.  condition. The latter implies the former, but the reverse is not true. 
Regarding the full rank a.e. kernel property, a close analog to Proposition \ref{theorem:rank_equivals}, although weaker for it provides only a sufficient condition, is the following result. 
\newtheorem{ram_th3}[rank_equivals]{Proposition}
\begin{ram_th3}\label{theorem:psisum_fullrank}
Let $\psi:U\times V \to \mathbb{R}$ admit the representation 
\begin{equation}\label{eq:psisum}
\psi(x,y) = \sum_{s=1}^{\infty} \phi_s(x)\xi_s(y), 
\end{equation}
for locally linear independent $\xi_s(y)$ such that for any $k\in\mathbb{N}$, there is a set of $k$ functions $\phi_s(x)$ that are locally linear independent in $U$. Then $\psi$ has full rank a.e. in $U\times V$. 
\end{ram_th3}

Proposition \ref{theorem:psisum_fullrank} states only a sufficient full rank condition. It turns out that for analytic in both arguments kernels, this condition becomes necessary also. In fact, the next result is a corollary of Propositions \ref{theorem:rank_equivals} and \ref{theorem:psisum_fullrank}.
\newtheorem{cor_analytic_full_rank_cond}{Corollary}
\begin{cor_analytic_full_rank_cond}\label{cor:analytic_full_rank_cond}
If $\psi$ is analytic in $U\times V$, then $\psi$ has full rank $a.e.$ in $U\times V$ if and only if 
for the Taylor functions $c_s(x)$ of $\psi$ about any point $p\in V$ and every $k\in\mathbb{N}$, there is a set of $k$ functions $c_s(x)$ that are locally linear independent in $U$. 
\end{cor_analytic_full_rank_cond}

We can apply the above result to a large class of symmetric analytic functions. Let $x.y = \sum_{l=1}^n x_{[l]} y_{[l]}$ denote the dot product in $\mathbb{R}^n$. Consider the class of functions $\psi(x,y) = h(x.y)$, $(x,y)\in U\times V$, where 
$W=\{x.y|x\in U, y\in V\}\subset \mathbb{R}$, $0\in V$, and 
$h$ is analytic function in $W$ such that $h^{(s)}(0)\ne 0$ for infinitely many $s$. 
The Taylor functions of $\psi(x,y)$ at $y=0$ are 
$$
c_s(x) = \frac{1}{[s]!} \frac{\partial^s}{\partial y^s} \psi(x,y)|_{y=0} = \frac{h^{(s)}(0)}{[s]!} x^s, 
$$
and since any finite collection of power functions is locally linear independent, 
for any $k$, there exist $k$ functions $c_s$ that are locally linear independent in $U$. 
Then, by Proposition \ref{theorem:psisum_fullrank}, $\psi$ has full rank $a.e.$ in $U \times V$. 
As an illustration, one can show the full rank a.e. property of the Laplace kernel $\psi(x,y) = \exp(-x.y)$. 
Other examples of kernels of the same kind are $cos(x.y)$ and $arccos(x.y)$. 

For subsets $U_l$, $l=1,...,m$, of $\mathbb{R}^n$, with $\otimes_{l=1}^{m} U_l$ we shall denote 
the product $U_1 \times ... \times U_m$. 
Let $\mu$ be the Lebesgue measure in $\mathbb{R}^n$. For the product $\otimes_{l=1}^m U_l$ of measurable $U_l$, by $\mu(\otimes_l U_l)$ we shall understand its Lebesgue (product) measure as a measurable subset of $\mathbb{R}^{nm}$.
To motivate the usage of ``a.e." notion in our definition of full rank, we need the following measure-theoretic result, which is of interest by itself. 
\newtheorem{ram_th6}[rank_equivals]{Proposition}
\begin{ram_th6}\label{theorema:f_nullset_zeros}
Let $f$ be an analytic function in an open subset $U$ of $\mathbb{R}^n$ 
such that the set $V_0 := \{x\in U| f(x) = 0\}$ has no interior points, then $V_0$ is measurable and $\mu(V_0) = 0$.
\end{ram_th6}

Note that if the function $\psi(x,y)$ is analytic in $U \times V$, then for any $k\ge 1$,  
$f(x_1,...,x_k,y_1,...,y_k) = det(\{\psi(x_i,y_j)\}_{i,j=1}^{k,k})$ is an analytic function 
in $(\otimes_{l=1}^{k} U \otimes_{l=1}^{k} V)$. 
If in addition, $\psi$ has a full rank $a.e.$ in $U\times V$, then, 
by the virtue of Proposition \ref{theorema:f_nullset_zeros}, we have immediately the following result.
\newtheorem{cor_full_rank_zero_measure}[cor_analytic_full_rank_cond]{Corollary}
\begin{cor_full_rank_zero_measure}
If $\psi: U\times V\to\mathbb{R}$ is analytic in $V$ and has full rank $a.e.$, then for any $k\ge 1$ the set 
$$
D_k := \{(x_1,...x_k,y_1,...,y_k)\in \otimes_{l=1}^{k} U \otimes_{l=1}^{k} V \textrm{, s.t. } rank(\{\psi(x_i,y_j)\}_{i,j=1}^{k,k}) < k\}
$$
has (product) measure zero, $\mu(D_k) = 0$.
\label{cor:full_rank_zero_measure}
\end{cor_full_rank_zero_measure}
Indeed, the set $D_k$ is closed and has no interior points, for otherwise $det(\{\psi(x_i,y_j)\}_{i,j=1}^{k})$ will vanish in an open subset and $\psi$ will not have full rank a.e.. 
Essentially, this corollary motivates the use of ``a.e.'' notion in the full rank definition.

Finally, we present a non-trivial example of full rank a.e. kernel - the standard distance on the unit sphere. 

\newtheorem{ram_th5}[rank_equivals]{Proposition}
\begin{ram_th5}\label{theorem:full_rank_sphere}
The standard distance on the $n$-sphere, $\mathbb{S}^n$, $n\ge 2$, given by 
$$
arccos(p.q), \\ p,q\in \mathbb{S}^n \subset \mathbb{R}^{n+1}
$$
and its square have full rank a.e. on $\mathbb{S}^n$. 
\end{ram_th5}

The above result has a probabilistic formulation as well: for every $k$, the square matrix with entries $arccos^2(p_i.q_j)$ of random points $p_i$ and $q_j$, $i,j=1,...,k$, on a unit sphere $\mathcal{S}^n$, $n\ge 2$, sampled independently by a continuous distribution on the sphere, is non-singular with probability one. This is evident by the completeness of the Lebesgue measure. More specifically, by continuous distribution we understand one that has a density, i.e. is absolute continuous with respect to the Lebesgue measure. Then, the probabilistic formulation above follows from Corollary \ref{cor:full_rank_zero_measure}.

\newtheorem{ram_remark2}{Remark}
\begin{ram_remark2}
The full rank a.e. property of $arccos(x.y), \textrm{ } x,y\in \mathbb{S}^n$ can be inferred from another standard result in the linear operator theory. In the Hilbert space $L_2(\mathbb{S}^n)$ of square integrable functions on $\mathbb{S}^n$, the operator with kernel $arccos(x.y)$ is symmetric and thus $arccos(x.y) = \sum_{k\ge 1}\lambda_k\phi_k(x)\phi_k(y)$, where $\lambda_k$ and $\phi_k$ are the eigenvalues and eigenvectors of this operator 
(see the Theorems of Hilbert and Schmidt, \cite{riesz-nagy}, Sec. 97). 
Since $\phi_k$'s form an orthonormal system in $L_2(\mathbb{S}^n)$, they are necessarily linear independent. 
Moreover, one can expect that they are locally linear independent, though one still has to show it. 
Our approach however, has an advantage in two aspects: (1) finding the Taylor functions of $arccos$ seems simpler than finding the spectral functions $\phi_k$, and (2), in general, we do not need the Hilbert space assumption. 
\end{ram_remark2}

\section{Concluding Remarks}

We return to the question stated in the introduction: when can the solution of an integral equation (\ref{eq:ram_fredholm}) be approximated by discrete solutions of the linear system of equations (\ref{eq:ram_fredholm_discrete})? The answer will be positive if, first, (\ref{eq:ram_fredholm}) has an unique solution $f$, and second, the series of solutions $\hat{f}_k$ of (\ref{eq:ram_fredholm_discrete}) converge to $f$ in an appropriate sense, a.s. or in some $L^p$-norm. 
The full rank a.e. condition on the kernel $\psi$ guarantees the existence of $\hat{f}_k$, but unfortunately, is not enough for the invertibility of the operator $A_{\psi}$, for equation (\ref{eq:ram_fredholm}) may still have multiple solutions. This is illustrated by the next example. Let for $(x,y)\in(-\infty,\infty)\times(0,\infty)$ 
$$
\psi(x,y) = 1 + \sum_{s=1}^{\infty} \frac{x^s}{s!} (\frac{y}{2s} - 1)y^{2s-1}.
$$
This kernel $\psi$ has the form (\ref{eq:psi_analytic}) with Taylor functions $c_{2s-1}(x) = -\frac{x^s}{n!}$ and $c_{2s}(x) = \frac{x^s}{2sn!}$. 
The conditions in Proposition \ref{theorem:psisum_fullrank} for $\psi$ are satisfied and thus, $\psi$ has full rank a.e.. On the other hand, $c_s(x)$ are not linearly independent and one easily finds that $A_{\psi}$ is not invertible. Indeed, for $f(y) = e^{-y}$, we have $(A_{\psi}f)(x) = 0$, for all $x\in\mathbb{R}$, and therefore, 0 is an eigenvalue of the operator $A_{\psi}$. 

The invertibility of $A_{\psi}$ however, can be guaranteed under some additional conditions. Let $\psi$ be analytic in $y\in V$ and the Taylor functions $c_s(x)$ in the analytical development (\ref{eq:psi_analytic}) of $\psi$ be linearly independent in $U$. Let $r>0$ be the convergence radius of the Taylor expansion of $\psi$ in $y$. 
For $r<\infty$, consider the vector space $\mathbb{P}(V)$ 
$$
\mathbb{P}(V) := \{\textrm{all measurable } f:V\to\mathbb{R}, \textrm{ s.t. } \sup_{s\ge 1} \{ r^{-s} \int_V |f(y)y^s| dy \} < \infty \}
$$
equipped with the $L^1$ norm and let $\mathbb{P}(V)=L^1(V)$, if $r=\infty$. 
Assume that $A_{\psi}$ is a bounded and closed operator from $\mathbb{P}(V)$ to $L^1(U)$. A sufficient condition for the latter is $\sup_{\rho\in(0,r)} \sum_{s=1}^{\infty} \rho^{s}\int_U |c_s(x)|dx < \infty$. Then, $A_{\psi}$ cannot have 0 as an eigenvalue because the equation 
$$
(A_{\psi}f)(x) = \sum_{s=1}^{\infty} c_s(x) \int_{V} y^sf(y)dy = 0 
$$
has no non-zero solution for $f$; otherwise $c_s$ will not be independent. 
In this case, $A_{\psi}: \mathbb{P}(V) \to Im(A_{\psi}) \subset L^1(U)$ is a one-to-one map. 
Consequently, for any $g\in Im(A_{\psi})$, $A_{\psi}f = g$ has a unique solution $f$. 
If in addition $\psi$ has full rank a.e., we can construct a series of the discrete invertible operators $A_k\in\mathbb{R}^{k^2}$, $k=1,2,...$, that converge uniformly to $A_{\psi}$, $||A_{k} - A_{\psi}|| \to 0$ (recall the setup in the introduction). Let also $g_k\in\mathbb{R}^k$, such that $||g_k-g||\to 0$ and $f_k := A_{k}^{-1}g_k\in\mathbb{R}^k$.  Now, for any convergent subsequence $f_{k'}$, $k'\subset k$, $||f_{k'}-f_0||\to 0$, by the uniform convergence of $A_k$, we have $||A_{k'}f_{k'} - A_{\psi}f_0|| \to 0$. On the other hand, $A_{k'}f_{k'} = g_{k'}$, and hence, $A_{\psi}f_0 = g$ and $f_0=f$. Therefore, $||f_{k'} - f||\to 0$ and we have the desired convergence to the unique solution of $A_{\psi}f = g$. 
Unfortunately, 
we cannot guarantee the existence of a convergent subsequence of $f_k$. In fact, if $A_{\psi}$ is compact, then there will be denumerable many eigenvalues concentrating at 0, $A_{\psi}^{-1}$ will be unbounded and the problem of approximating $f$ will be severely ill-posed - the ill-posedness of the Fredholm integral equations of first kind is a well known problem. Nevertheless, we find the principle possibility of inverting $A_{\psi}$ and all of its discrete approximations, the matrices $A_k$ obtained by random quadrature, to be sufficient motivation for this study  and deserving some further investigations. 

Note that if $\psi$ is analytic in $V$ and its Taylor functions $c_s(x)$ are locally linear independent in $U$, then $\psi$ is of full rank a.e. and any restriction $A_{\psi}|_{x\in W}$, for open $W\subset U$, is an injective operator - this is what we understand by local invertibility. With other words, any function $f\in\mathbb{P}(V)$ is possible to be recovered from $g=A_{\psi}f$, if its image $g$ is only known in some arbitrary open $W$. Even when is intractable with quadrature in practice, because of the ill-posedness, such recovery remains a principle possibility that eventually can be achieved, for example, with regularization methods targeting restricted classes of functions. As seen from the proof of Proposition \ref{theorem:full_rank_sphere}, the squared spherical distance is an analytical kernel with locally linear independent Taylor functions, and hence, is an example of full rank a.e. and locally invertible integral operator. 


\appendix

\section{Proofs}

\begin{proof}[Proof of Lemma \ref{lemma:linindep_functions}]
We show the claim by induction. Since $f_1(x)$ is not identically zero, there is $x_1\in U$ such that $f_1(x_1) \ne 0$. 
Let us assume that we have $k-1$ points $x_1$,...,$x_{k-1}$ in $U$, such that $det(\{f_i(x_j)\}_{i,j=1}^{k-1}) \ne 0$ 
and for $s=1,...,k$ define $A_s:=det(\{f_i(x_j)\}, i=1,...,k, i \ne s; j = 1,...,k-1)$. Since $A_{k-1} \ne 0$ and $f_i$ are linearly independent, we can choose $x_k\in U$  such that for $x=x_k$ 
$$
det \left( \begin{array}{cccc}
f_1(x_1) & f_2(x_1) & ... & f_{k}(x_1) \\
... & ... & ... & ... \\
f_1(x_{k-1}) & f_2(x_{k-1}) & ... & f_{k}(x_{k-1}) \\

f_1(x) & f_2(x) & ... & f_{k}(x) \\

\end{array} \right) = \sum_{i=1}^{k} (-1)^{i+k} A_{i} f_i(x) \ne 0,
$$
which proves the claim. 
\end{proof}
\begin{proof}[Proof of Proposition \ref{theorem:rank_equivals}]
We shall show that $(3) \Rightarrow (2) \Rightarrow (1) \Rightarrow (3)$.
Let condition $(3)$ hold and $\phi_1$, ..., $\phi_k$ be a basis of $span\{c_1, c_2,...\}$. Then, for any $s$, $c_s(x) = \sum_{i=1}^k \beta_s^i \phi_i(x)$ and $\psi(x,y) = \sum_{i=1}^k (\sum_s \beta_s^i (y-p)^s) \phi_i(x)$, provided that 
for any $i$ and $y\in V$, $\sum_s \beta_s^i (y-p)^s$ converges pointwise. 
Suppose the converse, that there exists $i_0$ and $y_0\in V $ such that $\sum_s \beta_s^{i_0} (y_0-p)^s$ does not converge.
Since $\phi_i$'s are linearly independent, we can choose $x_j\in U$, $j=1,...,k$, such that $rank(A := \{\phi_i(x_j)\}_{j,i=1}^{k}) = k$. Let $||A^{-1}|| > 0$ denote the operator norm of matrix $A^{-1}$. Fix a number $\epsilon > 0$. 
Since all $\sum_s c_s(x_j) (y_0-p)^s$ converge absolutly, there is $N$ such that 
$$
|\sum_{[s]=N}^M (\sum_{i=1}^k \beta_s^i \phi_i(x_j)) (y_0-p)^s| < \epsilon,
$$ 
for any $j$ and $M\ge N$. On the other hand, by the divergence assumption, there is $M > N$, such that 
$$
|\sum_{[s]=N}^M \beta_s^{i_0} (y_0-p)^s| \ge \sqrt{k} ||A^{-1}|| \epsilon.
$$ 
Define the k-vectors $z=(z_1,...,z_k)^T$ and $w=(w_1, ..., w_k)^T$, where 
$$
z_i := \sum_{[s]=N}^M \beta_s^{i} (y_0-p)^s
\textrm{ and } 
w_j := \sum_{[s]=N}^M (\sum_{i=1}^k \beta_s^{i} \phi_i(x_j)) (y_0-p)^s.
$$ 
The system $A z = w$ can be solved for $z$, $z = A^{-1}w$ and therefore $||z|| \le ||A^{-1}||\textrm{ }||w|| < \sqrt{k} ||A^{-1}|| \epsilon$, 
which contradicts $||z|| \ge |z_{i_0}| \ge \sqrt{k} ||A^{-1}|| \epsilon$.
Hence, the initial assumption is false and all 
$$
g_i(y) := \sum_{s=1}^{\infty} \beta_s^i (y-p)^s,
$$
are well defined functions in $V$ and condition (2) holds.

That (1) follows from (2) is immediate. Apparently, the rank of $\psi$ can not be larger than $k$. 
Moreover, by the linear independence assumption and Lemma \ref{lemma:linindep_functions}, there are $x_j\in U$ and $y_j\in V$, $j=1,...,k$, such that $rank(A := \{\phi_i(x_j)\}_{i,j=1}^k)=k$ and $rank(B = \{\xi_i(y_j)\}_{i,j=1}^k)=k$.
Therefore, for $\Psi:=\{\psi(x_i,y_j)\}_{i,j=1}^{k,k}$ we have $\Psi = AB^T$, which implies $rank(\Psi) = k$. 

Next we show $(1) \Rightarrow (3)$. Let $rank(\psi) = k$. For fixed $x_i \in U$, $i=1,...,k+1$, consider the set 
$\{\psi(x_1, y),...,\psi(x_{k+1}, y)\}$ of functions in $V$. 
Since for any $k+1$ points $y_j\in V$, no full rank matrix $\{\psi(x_i,y_j)\}_{i,j=1}^{k+1}$ exists, by Lemma \ref{lemma:linindep_functions}, $\{\psi(x_i, y)\}_{i=1}^{k+1}$ are not linear independent.
Then, there are functions $\alpha_i(\underline{x})$ of $\underline{x}=(x_1,...,x_{k+1})$, not all zero for any $\underline{x}$, such that 
$$
\sum_{i=1}^{k+1} \alpha_i(\underline{x}) \psi(x_i, y) = 0, \forall y\in V.
$$
By Taylor expanding $\psi$ in $y$ about a point $p\in V$ we obtain
$$
\sum_{s} \sum_{i=1}^{k+1} \alpha_i(\underline{x})c_s(x_i) (y-p)^s = 0, \forall y\in V,
$$
implying $\sum_{i=1}^{k+1} \alpha_i(\underline{x}) c_s(x_i) = 0, \forall s$.
The latter is true for any $x_i\in U$, $i=1,...,k+1$ and therefore, again by Lemma \ref{lemma:linindep_functions}, no set of $k+1$ functions $c_s$ is linearly independent. On the other hand, there is a set of $k$ functions $c_s$ that is linearly independent, for otherwise, following $(3)\Rightarrow (2) \Rightarrow (1)$, we would have that $rank(\psi) < k$. 
Consequently, there exists a set of functions $\phi_1(x)$,...,$\phi_k(x)$ among $c_s$, such that for any $s$, 
$c_s(x) = \sum_{i=1}^k \beta_s^i \phi_i(x)$, i.e. $c_s\in span\{\phi_1,...,\phi_k\}$. Thus, the condition (3) is fulfilled. 
\end{proof}
\begin{proof}[Proof of Lemma \ref{lemma:lli_powers}]
We will show the claim by induction on the number of dimensions. 
Let us first convince ourselves that it is true for $n=1$ and $U=(a-\delta,a+\delta)$, $\delta>0$. 
If $\sum_{l=1}^{\infty} \alpha_l x^{s_l} = 0$, $\forall x\in(a-\delta,a+\delta)$, 
then the above power series has radius of convergence at least $\rho:=\max\{|a-\delta|, |a+\delta|\}>0$. Hence, $f := \sum_{l=1}^{\infty} \alpha_l x^{s_l}$, being an analytic function in $(-\rho, \rho)$, vanishes everywhere in $(-\rho, \rho)$, which leads to the only possible choice $\alpha_l = 0$, for all $l$. 

Now, let $n>1$, $x=(x_{[1]},...,x_{[n]})$, and assume that all sets of power functions in $n-1$ dimensions are locally linear independent. If $\sum_l \alpha_l x^{s_l} = 0$ in $U\subset\mathbb{R}^n$, then the series will be absolute convergent in $U$ (we may eventually need to shrink $U$ a bit so that the closure $\bar{U}$ is entirely within the radius of convergence of the series). Consequently, we can regroup by $x_{[1]}$: $\sum_l \alpha_l x^{s_l} = \sum_{k=1}^{\infty} (\sum_{l: s_l^1=k} \alpha_l \underline{x}^{s_l^{n-1}} ) x_{[1]}^k$, where $\underline{x} := (x_{[2]},...,x_{[n]})$, $s_l=[s_l^1, s_l^{n-1}]$ and $s_l^{n-1}$ is a $(n-1)$-multi-index. Applying the argument from the first part of the proof, we have $\sum_{l: s_l^1=k} \alpha_l \underline{x}^{s_l^{n-1}} = 0$ in $U$ for every $k$. 
Finally, by the assumption, the only possible choice is $\alpha_l = 0$, for all $l$.
\end{proof}
\begin{proof}[Proof of Lemma \ref{lemma:lli_cond}]
We can repeat the lines of the proof of Lemma \ref{lemma:linindep_functions} with 
the only change at each selection step: choosing $x_i \in U_i$ instead of $x_i \in U$. 

Since $f_1(x)$ can not be identically zero in $U_1$, there is $x_1\in U_1$ such that $f_1(x_1) \ne 0$. 
Assuming that there are $x_i\in U_i$, $i=1,..,k-1$, such that $det(A_k := \{f_i(x_j)\}_{i,j=1}^{k-1}) \ne 0$, 
we can choose $x_k\in U_k$  such that for $x=x_k$ 
$$
det \left( \begin{array}{cccc}
f_1(x_1) & f_2(x_1) & ... & f_{k}(x_1) \\
... & ... & ... & ... \\
f_1(x_{k-1}) & f_2(x_{k-1}) & ... & f_{k}(x_{k-1}) \\
f_1(x) & f_2(x) & ... & f_{k}(x) \\
\end{array} \right) = \sum_{i=1}^{k} (-1)^{i+k} A_i f_i(x) \ne 0. 
$$
Otherwise, if such choice is impossible, $\sum_{i=1}^{k} A_i f_i(x) = 0$, $\forall x\in U_k$ and $f_i$ would not be locally linear independent.
\end{proof}
\begin{proof}[Proof of Proposition \ref{theorem:psisum_fullrank}]
The convergence of the right-hand series in (\ref{eq:psisum}) is assumed to be in absolute sense. Thus, $\sum_{s=1}^{\infty} |\phi_s(x)\xi_s(y)| < \infty$, for all $x\in U$ and $y\in V$. 

Suppose that for a number $k$ there exist open sets $U_i\subset U$ and $V_j\subset V$, $i,j=1,...,k$, 
such that for all $x_i\in U_i$ and $y_j\in V_j$, $rank(\{\psi(x_i,y_j)\}_{i,j=1}^{k,k}) < k$.
If we fix $x_i\in U_i$, $i=1,...,k$, then, by Lemma \ref{lemma:lli_cond}, 
$\{\psi(x_i, y)\}_{i=1}^{k}$ are not locally linear independent functions of y in $V$. 
Consequently, for the fixed $\underline{x}=(x_1,...,x_{k})$, there exists an open set $W(\underline{x})\subset V$, where $\{\psi(x_i, y)|_{W(\underline{x})}\}_{i=1}^{k}$ are linearly dependent, or equivalently, 
there are $\alpha_i(\underline{x})$, not all zero, 
such that 
$$
\sum_{i=1}^{k} \alpha_i(\underline{x}) \psi(x_i, y) = 0, \forall y\in W(\underline{x}).
$$
Taking into account the given representation (\ref{eq:psisum}) of $\psi$, we can write 
$$
\sum_{s=1}^{\infty} \sum_{i=1}^{k} \alpha_i(\underline{x}) \phi_s(x_i) \xi_s(y) = 0, \forall y\in W(\underline{x}).
$$
The local linear independence of $\xi_s$ then guarantees that  
$
\sum_{i=1}^{k} \alpha_i(\underline{x}) \phi_s(x_i) = 0, \forall s.
$
Therefore, for any choice of indices $s_1,...,s_k$, $rank(\{\phi_{s_j}(x_i)\}_{i,j=1}^{k}) < k$. 
Since the latter is true for any $x_i\in U_i$, $i=1,...,k$, 
by Lemma \ref{lemma:lli_cond}, no $k$ functions $\phi_s$ would be locally linear independent in $U$, 
a contradiction with the proposition assumption. 
Therefore, $\psi$ has full rank a.e. in $U\times V$. 
\end{proof}
\begin{proof}[Proof of Corollary \ref{cor:analytic_full_rank_cond}] 
By $\psi$ being analytic in $U\times V$ we mean that for every $p\in V$, $\psi(x,y) = \sum_s c_s(x) (y-p)^s$, for analytic in $U$ functions $c_s(x)$. The sufficiency then follows from Proposition \ref{theorem:psisum_fullrank}, 
which can be applied because $(y-p)^s$ are locally linear independent in $V$. 

Next, we shall check the necessity. Let $\psi$ have full rank $a.e.$ in $U\times V$ and there is $p\in V$ and $k$ such that there is no a collection of $k$ Taylor functions $c_s$ that is locally linear independent. 
But the local linear independence of the analytic functions $c_s$ is equivalent to their simple linear independence and thus, 
there is not set of $k$ functions $c_s$ that is linear independent in $U$. Then, $span\{c_1, c_2,...\}$ has a finite basis of size at most $k$ in $U$ and, by Proposition \ref{theorem:rank_equivals}, $rank(\psi) \le k$. 
This contradiction with the full rank assumption on $\psi$ proves the necessity. 
\end{proof}
For what follows, we need the following fact. 
\newtheorem{ram_lemma4}[lemma_linindep]{Lemma}
\begin{ram_lemma4}\label{lemma:union_open_sets}
Any collection $\mathcal{O}$ of (nonempty) open and disjoint subsets of $\mathbb{R}^n$ is countable.
\end{ram_lemma4}
\begin{proof}[Proof of Lemma \ref{lemma:union_open_sets}]
Let $O\in\mathcal{O}$. Then there is at least one vector $q\in O$ with all rational coordinates. 
By applying the axiom of choice we define a map 
$
\theta: O \mapsto (k_1,l_1,k_2,l_2,...,k_n,l_n),
$
where $k_i$, $l_i$, $i=1,...,n$, are integers such that 
$q=(2^{k_1}l_1,...,2^{k_n}l_n) \in O$. 
$\mathcal{O}$ is countable since $\theta$ is injective map 
from $\mathcal{O}$ to $\mathbb{Z}^{2n}$.
\end{proof}
\begin{proof}[Proof of Proposition \ref{theorema:f_nullset_zeros}]
Define $V_1 := \{x\in U| \frac{\partial f}{\partial x}|_{x} = 0\}$.
When $n>1$, $\frac{\partial f}{\partial x}$ is a gradient vector and $\frac{\partial f}{\partial x}|_{x} = 0$ 
means that all partial derivatives at point $x$ vanish. 
Both $V_0$ and $V_1$ are closed in $U$. We shall show that $V_0\backslash V_1$ has measure zero.

By the implicit function theorem, for any point $x\in V_0\backslash V_1$ there is 
an open ball $B_x=B_x(r_x)$ of radius $r_x>0$ and 
a hypersurface $S_x$ of dimension $n-1$, the graph of a function of $n-1$ variables, 
such that $x\in S_x = B_x \cap V_0\backslash V_1$. We have $\mu(S_x) = 0$.
Note that any two $S_x$ and $S_y$ that intersect and are images of open subsets of a common coordinate plane, say $\{x_{[i]}=0\}$, can be combined to the image of a continuous function in $\{x_{[i]}=0\}$. We have therefore the representation 
$V_0\backslash V_1 = \cup_{i=1}^n \cup_{\alpha\in A_i} S_{\alpha}^i$, where all $S_{\alpha}^i$ are disjoined images of continuous functions 
in $\{x_{[i]}=0\}$ and as such $\mu(S_{\alpha}^i) = 0$, $\forall i$ and $\forall \alpha \in A_i$. 

Let $U_{\alpha}^i = \cup_{x\in S_{\alpha}^i} B_x(r_x/2)$. Then it must be that $U_{\alpha}^i \cap U_{\beta}^i = \emptyset$, $\forall \alpha\ne\beta\in A_i$.  
For otherwise, there would be $x\in S_{\alpha}^i$ and $y\in S_{\beta}^i$ such that $B_x(r_x/2)\cap B_y(r_y/2) \ne \emptyset$ and  
consequently, either $x\in B_y(r_y)$ or $y\in B_x(r_x)$. If for example $y\in B_x(r_x)$, then $y\in S_x$ and $y \in S_x \cap S_y$,  and therefore $S_{\alpha}^i \cap S_{\beta}^i \ne \emptyset$, which would be a contradiction. 
Hence, by the virtue of Lemma \ref{lemma:union_open_sets}, for every $i$, the disjoint open sets $\{U_{\alpha}^i\}_{\alpha}$ are at most countably many, that is, the index sets $A_i$ are countable. 
Therefore, $\mu(V_0\backslash V_1) = 0$ and $\mu(V_0) = \mu(V_1\cap V_0)$. 

We can repeat the same analysis for higher derivatives of order $l > 1$ 
and apply induction on $l$. Define 
$$
V_l := \{x\in U \textrm{ s.t. }  \frac{\partial^s f}{\partial x^s}|_{x} = 0 \textrm{ for all }s=(s_1,...,s_n) \textrm{ for which } [s] = l\}.
$$
We claim that $\mu(V_{l-1}\backslash V_l) = 0$ for all $l>0$. We have already showed the claim for $l=1$.

For any $s=(s_1,...,s_n)\textrm{ for which }[s] = l-1$, define 
$$
V_{l-1}^s := \{x\in U \textrm{ s.t. } \frac{\partial^s f}{\partial x^s}|_{x} = 0\}, ~  
V_l^s := \{x\in U \textrm{ s.t. } \frac{\partial^{s+1} f}{\partial x^s\partial x^j}|_{x} = 0,\textrm{ for all }j=1,...,n\}
$$
and observe that $V_{l-1} = \cap_{s:[s]=l-1} V_{l-1}^s$ and $V_{l} = \cap_{s:[s]=l-1} V_{l}^s$. 
The argument from the first part of the proof for showing $\nu_n(V_0) = \nu_n(V_1\cap V_0)$ is applicable 
to all pairs $(V_{l-1}^s, V_l^s)$ and consequently 
$\mu(V_{l-1}^s) = \mu(V_l^s\cap V_{l-1}^s)$, $\forall s,~ [s] = l-1$. 
Next, we use the fact that for any four measurable sets $A$, $B$, $C$ and $D$, satisfying $B\subset A$, $\mu(B)=\mu(A)$, 
$D\subset C$ and $\mu(D)=\mu(C)$, we have $\mu(B\cap D)=\mu(A\cap C)$, 
to conclude that 
$\mu(\cap_{[s]=l-1} V_{l-1}^s) = \mu(\cap_{[s]=l-1} (V_{l-1}^s \cap V_l^s))$. 
Hence we have $\mu(V_{l-1}) = \mu(V_{l-1} \cap V_l)$ and $\mu(V_{l-1}\backslash V_l) = 0$.
Now, assuming that $\mu(V_0) = \mu(\cap_{m=0}^{l-1} V_m)$, from 
$$
\mu(\cap_{m=0}^{l-1} V_m) = \mu((\cap_{m=0}^{l-1} V_m) \backslash V_l) + \mu(\cap_{m=0}^{l} V_m) = \mu(\cap_{m=0}^{l} V_m) 
$$ 
we obtain $\mu(V_0) = \mu(\cap_{m=0}^{l} V_m)$ and by the induction principle the latter is true for every $l$.

Note that $\mu(\cap_{m=0}^{l} V_m) \searrow \mu(\cap_{m=0}^{\infty} V_m)$, as $l\to\infty$. 
Finally, we realize that the set $\cap_{m=0}^{\infty} V_m$ should be empty, 
for otherwise, with all vanishing derivatives at one point, 
the analytic function $f$ would vanish in an open subset of $U$, a contradiction.
Thus we conclude that $\mu(V_0) = \mu(\cap_{m=0}^{\infty} V_m) = 0$.
\end{proof}
\begin{proof}[Proof of Proposition \ref{theorem:full_rank_sphere}]
Consider the following parametrization of the north hemisphere about the north pole: 
$$
\phi: x = (x_1,x_2,...,x_n) \in \mathbb{O}^n \to p=(x_1, x_2, ..., x_n, (1-\sum_{j=1}^n x_j^2)^{1/2}) \in \mathbb{R}^{n+1}, 
$$
where $\mathbb{O}^n = \{x\in[0,1]^n \textrm{, } \sum_{i=1}^n x_i^2 \le 1$ is the closed $n$-ball. Then the inner product of $p=\phi(x)$ and $q=\phi(y)$ is 
$$
p.q = \phi(x).\phi(y) = \sum_{i=1}^n x_iy_i + (1-\sum_{j=1}^n x_j^2)^{1/2}(1-\sum_{j=1}^n y_j^2)^{1/2}.
$$
The proof of the desired claim is based on application of Proposition \ref{theorem:psisum_fullrank}. For the purpose we need to derive the Taylor functions $f_s(x)$ in the expansion of $arccos(p.q)$ in $y$ at $0\in\mathbb{R}^n$. Specifically, we need to show that for any $k$, there are $k$ locally linear independent $f_s(x)$ in $\mathbb{O}^n$. 
We will first show this property for the Taylor functions of 
$$
\psi(x,y) = arccos(\phi(x).\phi(y)) \textrm{, } x,y\in \mathbb{O}^n,
$$
and then for the square spherical distance $\psi^2$. Note that $\psi$ is analytic in $\mathbb{O}^n \times \mathbb{O}^n$.

The first simplification we are going to make is to fix $y_2 =...=y_n=0$ and consider $\psi(x,y)$ as a function of $y_1\in[0,1]$. 
Then, it suffices to show the local linear independence condition only for the partial derivatives $\frac{\partial^s \psi}{\partial y_1^s}$. 

For brevity, denote $\underline{x} = (\sum_{i=1}^n x_i^2)^{1/2}$ and $\underline{y} = (\sum_{i=1}^n y_i^2)^{1/2}$.
The form of the first few derivatives of $z=\phi(x).\phi(y)$ as a function of $y_1$  
$$
\frac{\partial z}{\partial y_1} = x_1 - \frac{(1-\underline{x}^2)^{1/2}}{(1-\underline{y}^2)^{1/2}} y_1 \textrm{, }
\frac{\partial^2 z}{\partial y_1^2} = -\frac{(1-\underline{x}^2)^{1/2}}{(1-\underline{y}^2)^{1/2}} - \frac{(1-\underline{x}^2)^{1/2}}{(1-\underline{y}^2)^{3/2}} y_1^2 \textrm{, }
$$
can be generalized (the proof is by induction) to the following equations
\begin{equation}\label{eq:z_partials}
\frac{\partial z}{\partial y_1}|_{y=0} = x_1 \textrm{, }
\frac{\partial^{2k+1} z}{\partial y_1^{2k+1}}|_{y=0} = 0 \textrm{, }
\frac{\partial^{2k} z}{\partial y_1^{2k}}|_{y=0} = -a_k (1-\underline{x}^2)^{1/2} \textrm{, } k\ge 1,
\end{equation}
for some integers $a_k > 0$.

On the other hand, the general form of the derivatives of $\psi(z) = arccos(z)$ 
$$
\frac{\partial^{2s} \psi}{\partial z^{2s}} = 
- \sum_{l=1}^{s} b_{2s,l} \frac{z^{2l-1}}{(1-z^2)^{(2s+2l-1)/2}} \textrm{, }
$$
$$
\frac{\partial^{2s+1} \psi}{\partial z^{2s+1}} = 
- \sum_{l=0}^{s} b_{2s+1,l} \frac{z^{2l}}{(1-z^2)^{(2s+2l+1)/2}},
$$
calculated at $z=(1-\underline{x}^2)^{1/2}$ (corresponding to $y=0$), can be summarized to 
\begin{equation}\label{eq:psi_partials}
\frac{\partial^s \psi}{\partial z^s}|_{z=(1-\underline{x}^2)^{1/2}} =  - \frac{1}{\underline{x}^s} 
\sum_{l=0}^{s-1} b_{s,l} \Big( \frac{(1-\underline{x}^2)^{1/2}}{\underline{x}} \Big)^l \textrm{, } s\ge 1,
\end{equation}
for some integers numbers $b_{s,l} \ge 0$ with $b_{s,s-1} > 0$. Formally, (\ref{eq:psi_partials}) can be verified by induction.

Next, we combine equations (\ref{eq:z_partials}) and (\ref{eq:psi_partials}) to express the partial derivatives of $\psi$ in $y_1$ 
\begin{equation}\label{eq:psi_partials_y1}
\frac{\partial^s \psi}{\partial y_1^s} =  \sum_{k=1}^s \frac{\partial^{k} \psi}{\partial z^{k}} 
\sum_{l_1+...+l_k = s} d_{l_1l_2...l_k} 
\frac{\partial^{l_1} z}{\partial y_1^{l_1}}...\frac{\partial^{l_k} z}{\partial y_1^{l_k}}, 
\end{equation}
where $d_{l_1l_2...l_k}$ are some integers and the second sum is over all $l_1\le...\le l_k$. For example 
$$
\frac{\partial \psi}{\partial y_1}|_{y=0} = \frac{\partial \psi}{\partial z}\frac{\partial z}{\partial y_1}  = -\frac{x_1}{\underline{x}} \textrm{, }
$$
$$
\frac{\partial^2 \psi}{\partial y_1^2}|_{y=0} = \frac{\partial^2 \psi}{\partial z^2} \Big(\frac{\partial z}{\partial y_1}\Big)^2  + 
\frac{\partial \psi}{\partial z}\frac{\partial^2 z}{\partial y_1^2} 
= -\frac{(1-\underline{x}^2)^{1/2}}{\underline{x}} ((\frac{x_1}{\underline{x}})^2 - 1) \textrm{, etc.}
$$
The complex form of these functions, which mix $x_1$ with $\underline{x}$, renders the verification of their local linear independence untractable. 
Fortunately, at this point we can do further simplifications by taking into account that $n>1$. If we fix $x_1=0$ and show that for any $k$, there are $k$ local linear independent functions $\frac{\partial^s \psi}{\partial y_1^s}|_{x_1=0, y=0} (x_2,...,x_n)$, our claim will be proven. Note that this step can not be taken if $n=1$. Indeed, we have already verified that the claim is not true for the 1-dimensional circle (recall the example at the end of Section \ref{sec:finite_rank}).

First, observe that $\frac{\partial^{2s+1} \psi}{\partial y_1^{2s+1}}|_{x_1=0, y=0} (x_2,...,x_n) = 0$, 
because $\frac{\partial^{2k+1} z}{\partial y_1^{2k+1}}|_{x_1=0,y=0} = 0$, for all $k\ge 0$, and all terms in Eq. (\ref{eq:psi_partials_y1}) must have at least one odd $l_j$. 
Second, plugging (\ref{eq:z_partials}) and (\ref{eq:psi_partials}) in (\ref{eq:psi_partials_y1}) and denoting 
$\tilde{d}_{l_1l_2...l_k} = a_{l_1} a_{l_2} ... a_{l_k} d_{(2l_1)(2l_2)...(2l_k)} $, we obtain 
$$
\frac{\partial^{2s} \psi}{\partial y_1^{2s}}|_{x_1=0, y=0} (x_2,...,x_n) = 
$$
$$
\sum_{k=1}^{2s} \Big[ \sum_{l=0}^{k-1} b_{k,l} \Big( \frac{(1-\underline{x}^2)^{1/2}}{\underline{x}} \Big)^l
\sum_{l_1+...+l_k = s} \tilde{d}_{l_1l_2...l_k} 
\frac{1}{\underline{x}^k} 
(1-\underline{x}^2)^{k/2} \Big] = 
P_{s}\Big( \frac{(1-\underline{x}^2)^{1/2}}{\underline{x}} \Big),
$$
where $P_{s}$ is a polynomial with integer coefficients of degree $4s-1$. Hence, by introducing the functions $h(x_2,...,x_n) = \frac{(1-\underline{x}^2)^{1/2}}{\underline{x}}$, we conclude that 
\begin{equation}\label{eq:psi_partials_x1y1}
\frac{\partial^{2s} \psi}{\partial y_1^{2s}}|_{x_1=0, y=0} (x_2,...,x_n) \in 
span\{h,...,h^{4s-1}\} \setminus span\{h,...,h^{4s-2}\} \textrm{, } s > 1.
\end{equation}
Finally, by Lemmas \ref{lemma:lli_bijection} and \ref{lemma:lli_lin_comb}, we assert that indeed, for every finite set of distinct even numbers $s$, the corresponding partial derivatives (\ref{eq:psi_partials_x1y1}) are locally linear independent. 
The conditions of Proposition \ref{theorem:psisum_fullrank} therefore hold and the spherical distance has full rank a.e. (in the north hemisphere).

Now we consider the square spherical distance. From the equality
$$
\frac{\partial^{s} \psi^2}{\partial y_1^{s}} = 2 \psi \frac{\partial^{s} \psi}{\partial y_1^{s}} + 
\sum_{k=1}^{s-1} e_{k}^s \frac{\partial^{k} \psi}{\partial y_1^{k}} \frac{\partial^{s-k} \psi}{\partial y_1^{s-k}},
$$
for some integers $e_k^s>0$, we conclude that at the point $(x_1=0, y=0)$
$$
\frac{\partial^{2s+1} \psi}{\partial y_1^{2s+1}} (x_2,...,x_n) = 0 \textrm{, }
$$
\begin{equation}\label{eq:psi2_partials_even}
\frac{\partial^{2s} \psi}{\partial y_1^{2s}} (x_2,...,x_n) = 
\psi((1-\underline{x}^2)^{1/2}) P_{s}\Big( \frac{(1-\underline{x}^2)^{1/2}}{\underline{x}} \Big) + Q_{s}\Big( \frac{(1-\underline{x}^2)^{1/2}}{\underline{x}} \Big),
\end{equation}
where $Q_{s}=\sum_{k=1}^{s-1} e_k^s P_k P_{s-k}$ is a polynomial with integer coefficients of degree $4s-2$. 

Let us assume that for some $k$, the first $k$ even derivatives (\ref{eq:psi2_partials_even}) are not locally linear independent. 
Let $t=\frac{(1-\underline{x}^2)^{1/2}}{\underline{x}}$. Then, there are some $c_s$ with $c_k\ne 0$, such that 
$$
\sum_{s=1}^k c_s \Big( P_s(t) arccos(\frac{t}{(1+t^2)^{1/2}}) + Q_{s}(t) \Big) = 0 \textrm{, } t\in W, 
$$
in an open (non-empty) $W\subset[0,1]$. This is equivalent to 
$$
arccos(\frac{t}{(1+t^2)^{1/2}}) = - \frac{Q(t)}{P(t)} \textrm{, } t\in W, 
$$
where $P=\sum_{s=1}^k c_s P_s$ and $Q=\sum_{s=1}^k c_s Q_s$ are polynomial of degrees $4k-1$ and $4k-2$, respectively. After taking the derivatives on both sides of the above equation and applying some algebra we obtain the equality
$(P(t))^2=(1+t^2)(P(t)Q'(t)-Q(t)P'(t))$. The latter however is impossible because of the degree disparity, 
$8k-2$ for $P^2$ and $8k-3$ for the left-hand side polynomial 
(we use that $degree(P(t)Q'(t)-Q(t)P'(t)) \le degree(P)+degree(Q)-2$).
This proves the local linear independence of the first $k$ even derivatives (\ref{eq:psi2_partials_even}).
Due to Proposition \ref{theorem:psisum_fullrank}, the square spherical distance has full rank a.e. in $\mathbb{S}^n$.
\end{proof}


\section{Bibliography}

\begin{thebibliography}{10}

\bibitem{riesz-nagy}
Riesz, F., and B. Sz.-Nagy, 
\newblock{Functional Analysis},
\newblock{1955, Ungar, New York}.


\end{thebibliography}
\end{document}